\documentclass[12pt,a4paper,intlimits]{amsart}

\usepackage{amssymb,latexsym}

\newcommand{\field}[1]{\mathbb{#1}}
\newcommand{\R}{\field{R}}
\newcommand{\C}{\field{C}}
\newcommand{\N}{\field{N}}

\newcommand{\eps}{\varepsilon}
\newcommand{\supp}{\mathrm{supp}\,}

\newtheorem{thm}{Theorem}[section]
\newtheorem{prop}[thm]{Proposition}

\newtheorem{lemma}[thm]{Lemma}

\begin{document}

\bibliographystyle{amsplain}

\title[Estimates for maximal operators]{Damped oscillatory integrals and
boundedness of maximal operators associated to mixed homogeneous
hypersurfaces}

\author[I. A. Ikromov]{Isroil A.Ikromov}
\address{Department of Mathematics, Samarkand State University, University
Boulevard 15, 703004, Samarkand, Uzbekistan}
\email{{\tt ikromov1@rambler.ru}}

\author[M. Kempe]{Michael Kempe}
\address{Mathematisches Seminar, C.A.-Universit\"at, Ludewig-Meyn-Stra\ss e
4, D-24098 Kiel, Germany}
\email{{\tt kempe@math.uni-kiel.de}}

\author[D. M\"uller]{Detlef M\"uller}
\address{Mathematisches Seminar, C.A.-Universit\"at, Ludewig-Meyn-Stra\ss e
4, D-24098 Kiel, Germany}
\email{{\tt mueller@math.uni-kiel.de}}

\thanks{We acknowledge the support for this work by the Deutsche 
Forschungsgemeinschaft and the hospitality of the Erwin Schr\"odinger Institute 
in Vienna, which the second and third author could enjoy during their stay at 
the program ''Combinatorial and Number-Theoretic Methods in Harmonic Analysis''.}

\begin{abstract} We study the boundedness problem  for maximal operators
in 3-dimen\-sional Euclidean space associated to hypersurfaces given as the 
graph of $c+f$, where $f$ is a mixed homogeneous function which is smooth 
away from the origin and $c$ is a constant.
 Assuming that the Gaussian curvature of this surface nowhere vanishes of 
 infinite order, we prove that the associated maximal operator is bounded on 
 $L^p(\R^3)$  whenever $p>h\ge 2$. Here $h$ denotes a
``height'' of the function $f$ defined in terms of its  maximum order of 
vanishing and the weights of homogeneity. This result generalizes 
a corresponding theorem on mixed homogeneous polynomial functions by 
A. Iosevich  and E. Sawyer. In particular, it shows that a certain 
``ellipticity'' conditon used by these authors is not necessary. If $c\neq 0$, 
our result is sharp.
\end{abstract}
\maketitle

\thispagestyle{empty}

\section{Introduction}

A by now ``classical'' theorem of real analysis is E.M. Stein`s maximal theorem
for spherical means on Euclidean space $\R^n$ . Stein's results in \cite{MR54:8133a} covered
the case $n\geq 3$, and the remaining $2$-dimensional case was later dealt
 with 
by J. Bourgain \cite{MR88f:42036}. These results became the starting point for the study of
various classes of maximal operators associated to subvarieties, such as
maximal operators defined by

\begin{equation}\label{defmaxop} Mg(x)=\sup_{t>0}\Bigl\vert\int_S
g(x-ty)\psi(y)d\sigma(y)\Bigr\vert,
\end{equation}
 where $S$ is a smooth hypersurface, $\psi$ is a fixed
non-negative function in $C_0^\infty(S)$, and $d\sigma$ the surface 
measure on $S$. For instance, A. Greenleaf \cite{MR84i:42030} proved that $M$ is
bounded on $L^p(\R^n)$, if $n\geq 3$ and $p>\frac{n}{n-1}$, provided $S$ has
everywhere non-vanishing Gaussian curvature and is star-shaped with respect
to the origin. In contrast, the case where the Gaussian curvature vanishes at some 
points is still widely open, and sharp results for this case are known only for
 particular classes of surfaces. A result of general nature is given by  
 C.D. Sogge and E.M. Stein in \cite{MR87d:42030}, where they show that if the Gaussian
curvature of $S$ does not vanish of infinite order at any point of $S$ then
$M$ is bounded on $L^p$ in a certain range $p>p(S)$. However, the exponent $p(S)$ in 
that paper is in general far from being optimal.

 It is well-known that the  $L^p$-estimates of the maximal operator (\ref{defmaxop})
 are strongly related 
to the decay of the Fourier transform of measures carried on $S$, i.e.\ to oscillatory 
integrals of the form
\begin{equation}\label{ftrafosurf}
\int_S e^{-i\xi \cdot x}\psi(x)\, d\sigma(x),
\end{equation} where $\psi$ is a compactly supported density on $S$. 

The decay of the oscillatory integral \eqref{ftrafosurf} as $|\xi|\to\infty$ in return 
is connected to geometric properties of the surface $S$ and has been
considered by various authors, including van der Corput, E. Hlawka, C.S.
Herz, W. Littman, B. Randol, I. Svensson, A. Varchenko, C.D. Sogge, E.M.
Stein, J.J. Duistermaat, Colin de Verdier,  et al.. We refer to \cite{MR95c:42002} for references, also 
to results on maximal operators associated to surfaces.

Another important idea, introduced in \cite{MR87d:42030} and applied in several subsequent 
articles,  is to ``damp'' the 
oscillatory integral \eqref{ftrafosurf}, by multiplying the amplitude $a$ by a 
 suitable power of the Gaussian curvature on $S$, in order to obtain the ``optimal'' 
decay of order $|\xi|^{-(n-1)/2}$.

A case which has been studied quite comprehensively is the case of convex
hypersurfaces of finite type. Sharp estimates for the Fourier tranform of
measures carried on such surfaces $S$ have been obtained by J. Bruna, A. Nagel
and S. Wainger \cite{MR89d:42023}. In that article, the authors introduce a family of 
nonisotropic balls on $S$, called ``caps'', by setting
\begin{equation*}
B(x,\delta):=\{y\in S\vert\, \mathrm{dist}(y,T_x S)<\delta\},
\end{equation*}
where $T_x S$ denotes the tangent space of $S$ at $x$.
Suppose that $\xi$ is normal to $S$ at $x_0$ and that the density $\psi$ is smooth with sufficiently small support.
 Then it is shown that the estimate

\begin{equation}\label{caps}
 |\widehat{\psi d\sigma}(\xi)|\leq C|B(x_0,|\xi|^{-1})|
\end{equation}
 holds, where $|B(x_0,\delta)|$ denotes the surface area of
$B(x_0,\delta)$.

The results by A. Nagel, A. Seeger and S. Wainger in \cite{MR94m:42033} on maximal
operators associated to such convex surfaces are based on this type of
estimates for Fourier transforms of surface-carried measures. Moreover they
obtain sharp results for convex hypersurfaces given as the graph of a
mixed homogeneous convex function $Q$ (i.e.\ there exist even integers
$(a_1,\dots,a_n)$ such that
$Q(s^\frac{1}{a_1}x_1,\dots,s^\frac{1}{a_n}x_n)=sQ(x),s>0)$, generalizing
the results by M. Cowling and G. Mauceri in \cite{MR87m:42013}.

Further results on the boundedness problem for maximal operators associated
to convex hypersurfaces were based on a result due to H.~Schulz \cite{MR93a:42007}
(see also  \cite{MR58:16649}), which states that, possibly after rotating
 the coordinate system, any smooth convex finite type function $\Phi$ can be
 written in the form $\Phi(x)=Q(x)+R(x)$, where $Q$ is a convex mixed homogeneous polynomial
that vanishes only at the origin, and $R(x)$ is a remainder term, in the
sense that it tends to zero under non-isotropic dilations of $Q$.

By using this result, A. Iosevich and E. Sawyer  \cite{MR99b:42023} proved that if
$p>2$ then the operator (1) is $L^p$ bounded if and only if the inequality
$p>(\frac{1}{a_1}+\dots +\frac{1}{a_n})^{-1}$ holds.

 However, it seems that the damping by powers of the Gaussian curvature does not
give the best possible results for the associated maximal operators, if the hypersurface 
is non-convex.

For this reason, M. Cowling and G. Mauceri \cite{MR89e:42013} and also A. Iosevich
and E. Sawyer \cite{MR99b:42023} have considered oscillatory integrals with another
damping factor, namely powers of the original phase function, if one expresses
\eqref{ftrafosurf} in terms of local coordinates for $S$.
In this way, they obtained sharp results for finite type convex
hypersurfaces.

In general, the boundedness problem for maximal operators associated to
non-convex hypersurfaces  is widely open.

This is partly due to the fact that, for non-convex hypersurfaces, Bruna, Nagel,
Wainger--type estimates (even Randol type estimates \cite{MR40:4678b}) like \eqref{caps}
fail to be true, even for non-convex smooth hypersurfaces in $\R^3$
(see \cite{MR95b:42021}).
\smallskip

In order to describe our results, let us introduce some notation.
 For any pair $\kappa=(\kappa_1,\kappa_2)$ of ''weights'' $\kappa_1,\kappa_2>0,$ 
 we define a
 group of dilations $\{\delta_r\}_{r>0}$ on $\R^2$ by setting
$\delta_r(x):=(r^{\kappa_1}x_1,r^{\kappa_2}x_2)$. A function $f$ on
$\R^2\setminus\{0\}$ is called $\kappa$-{\it homogeneous} of degree 
$\alpha\in \R,$
 if
\begin{equation*} f(\delta_r(x))=r^\alpha f(x)\quad\quad\text{for all}\  x\neq 0,
r>0.
\end{equation*} If $f$ is twice differentiable at $x,$ denote by $D^2 f(x)$
the Hessian matrix of $f$ at $x,$ and its determinant by
$\mathrm{Hess}f(x):=\det D^2f(x)$.

The {\it order} $\mathrm{ord} f(x)$ of $f$ in $x$ is understood to be the smallest
non-negative integer $j$ such that $D^j f(x)\neq 0,$ where $D^j f(x)$ denotes the
$j$-th order total derivative of $f$ in $x$. If $f$ is $\kappa$-homogeneous, then clearly

\begin{equation*}
\mathrm{ord}f:=\sup\limits_{x\in\R^2\setminus\{0\}}\mathrm{ord}
f(x)=\sup\limits_{x\in S^1}\mathrm{ord} f(x),
\end{equation*}
where $S^1$ denotes the unit circle in $\R^2.$

The {\it height} of a $\kappa$-homogeneous function $f$ is defined by

\begin{equation*}
h:=\max\Bigl\{\frac{1}{\kappa_1+\kappa_2},\mathrm{ord}f\Bigr\}.
\end{equation*}
We consider surfaces of the form 
$$S=\{(x_1,x_2,c+f(x_1,x_2))\vert\, (x_1,x_2)\in\R^2\}\subset\R^3,$$ 
where $f$ is a $\kappa$-homogeneous function of degree $1$ which is smooth 
away from the origin, and where $c\in\R$ is fixed. Let $M$ denote the corresponding 
maximal operator, i.e. 

\begin{equation*}
Mg(x):=\sup_{t>0}\Bigl\vert\int_Sg(x-ty)\psi(y)\, d\sigma(y)\Bigr\vert,\quad\quad
g\in C_0^\infty(\R^3),
\end{equation*} where $\psi\in C_0^\infty(S)$  and $d\sigma$ denotes the surface
measure on $S$.

Partial derivatives of differentiable functions are written in the form
$\partial^k_{j_1,\dots,j_k}f(x):=\frac{\partial^k f}{\partial x_{j_1}\dots\partial x_{j_k}}(x)$.

\begin{thm}\label{mainthm}
Assume that $f$ is $\kappa$-homogeneous of degree 
$1$ and smooth away from the origin, where $\kappa\neq(1,1)$, and that
$\mathrm{ord}(\mathrm{Hess}f)<\infty$ and $h\geq 2.$
Then the maximal operator $M$ is bounded on $L^p(\R^3)$ whenever $p>h$. 

Moreover, if $\psi(0)c\neq 0$, then the maximal operator is unbounded for $p\leq h$.
\end{thm}

{\bf Remarks.} a) This theorem is in fact a generalization of a theorem by 
A. Iosevich and E. Sawyer \cite{MR97f:42035}, who study the case of $\kappa$-homogeneous polynomial phase 
functions $f$ with isolated critical point at the origin such that $\mathrm{Hess}f$ 
and $f$ have no common zero, except for the origin.

b) We shall use damping of the arizing oscillatory integrals \eqref{ftrafosurf} by 
powers of a ``homogenized'' version of $|\nabla f(x)|$ instead of powers of $|f(x)|$.
 This has some technical advantages and seems, in a way, to be even more natural.

c) If $f$ is a $(1,1)$-homogeneous function of
degree $1$ then the surface $S$ is  conic.
In this case damped oscillatory integrals can not decay faster than
$O(|\xi|^{-\frac{1}{2}}),$ so that the boundedness
problem for associated maximal operators becomes much more complicated.
Nevertheless, if at every point one principal curvature of the conic surface does 
not vanish then the $L^p$-boundedness of the maximal operator for $p>2$ follows 
from C. Sogge's results in \cite{MR96e:42014}. In view of this, Theorem \ref{mainthm}
gives an almost complete answer to the question of $L^p$-boundedness of the maximal 
operator $M$, in the case where $h\ge 2.$ In the case where $h<2,$ damping of the 
oscillatory seems to be of no use, and new ideas seem to be needed in order to 
obtain sharp results. This case remains open at this time.
\medskip

The article is organized as follows: In Section \ref{reduction}, we outline some well-known
 reductions to oscillatory integral estimates. Basic for these estimates will be 
some structural results on $\kappa$-homogeneous functions, which will be derived in 
Section \ref{auxiliary}. Section \ref{oscestimates} will then contain the required estimates 
of the damped oscillatory intergrals which arize in the context of Theorem
\ref{mainthm}. The  proof of the sharpness of our main theorem will be carried out in Section
\ref{sharpness}.

\section{Preliminary reductions}\label{reduction}

Following a standard approach (see e.g. \cite{MR95c:42002}), the theorem will be
shown by embedding $M$ respectively linearizations of $M$  into an analytic
 family of operators and then
interpolating an $L^2\to L^2$ and an $L^\infty\to L^\infty$-estimate.

Define $\widetilde{\psi}\in C_0^\infty(\R^2)$ by $\widetilde{\psi}(x):=\psi(x,c+f(x))$.
First, we restrict ourselves to the situation where $\kappa_1\neq 1\neq\kappa_2$. The remaining
cases will be discussed afterwards. To start with,
we observe that the zeros of $\nabla f\vert_{S^1}$ form a discrete, hence finite subset 
of the unit circle. This follows immediately from Lemma \ref{localform}. Choose
a partition of unity  $\varphi_j$ on $S^1$,
$0\leq j\leq k$   such that at exactly one zero $x_j$ of
$\nabla f$ lies in $\supp\varphi_j$ if $1\leq j\leq k$, and $\nabla f$ does not vanish on $\supp\varphi_0$.

Obviously $\varphi_j$ admits exactly one continuation to a $\kappa$-homogeneous function of
degree $0$ on $\R^2\setminus\{0\}$, which is again denoted by $\varphi_j$. If we define
maximal operators $M_jg(x):=\sup_{t>0}\vert\int g(x-ty)\psi(y)\varphi_j(y)\,
d\sigma(y)\vert,$
then

\begin{equation*}
Mg(x)\leq\sum_{j=0}^k M_jg(x).
\end{equation*}
We can therefore assume without loss of generality that either $\nabla f$ does not vanish on $\supp\widetilde{\psi}$
or that there is exactly one $x^0\in S^1$ such that every $x\in\supp\widetilde{\psi}$ with $\nabla f(x)=0$ lies on the
$\kappa$-homogeneous curve $\{\delta_r(x^0)\vert\, r>0\}$ generated by $x^0$. Denote by $n$ the order of $f$
at that point and by $G_f\colon\R^2\setminus\{0\}\to\R$ the unique $\kappa$-homogeneous function of degree one
with $G_f(x)=\vert\nabla f(x)\vert^{\frac{n}{n-1}}$, $x\in S^1$.

For $w\in\C$ with $\mathrm{Re}\,w>-\frac{1}{h}$ define a measure $d\sigma_w(x):=G_f(x)^w
\psi(x)d\sigma(x)$ on $S$ as well as the corresponding maximal operator
 
\begin{equation*} M_w g(x):=\sup_{t>0}\Bigl\vert\int_S g(x-ty)\,
d\sigma_w(y)\Bigr\vert.
\end{equation*}
If $\nabla f$ does not vanish on $\supp\widetilde{\psi}$ we do not need any damping factor and just take $G_f\equiv 1$.

Obviously, $M_0=M$. It follows easily from the $\kappa$-homogeneity of $f$ and 
 Lemma \ref{localform} that 
$G_f^w$ is locally integrable if $\mathrm{Re\,}w>-\frac{1}{h}$, and thus $M_w$
is bounded on $L^\infty(\R^3)$ for these values of $w$.
Moreover, if we assume without loss of generality that $\psi\geq 0$, then
 $M_w g(x)\leq M_{\mathrm{Re}w} \vert g\vert(x)$. Once we can show
the $L^2$-boundedness of $M_\alpha$ for
$\alpha>\frac{1}{2}-\frac{1}{h}$, then Theorem \ref{mainthm} follows from 
Stein's interpolation theorem for analytic families of operators \cite{MR46:4102}.
But, it is an easy consequence of Sobolev's embedding theorem that 
 $M_\alpha$ is bounded on $L^2$
if the Fourier transform of the measure $d\sigma_\alpha$ satisfies

\begin{equation}\label{ft1}
\vert\widehat{d\sigma_\alpha}(\xi)\vert\leq C(1+\vert\xi\vert)^{-(\frac{1}{2}+
\eps)},
\end{equation} 
\begin{equation}\label{ft2}
\vert\nabla\widehat{d\sigma_\alpha}(\xi)\vert
\leq C(1+\vert\xi\vert)^{-(\frac{1}{2}+\eps)}
\end{equation}
for some $\eps>0$ (see e.g. \cite{MR87d:42030}). These estimates will be 
established in Theorem \ref{mainprop}.

If now for example $\kappa_1=1$ (the case $\kappa_2=1$ can be treated analogously), in view 
of Lemma \ref{localform} we can choose 
a more refined partition of unity in order to ensure that either $\partial_2 f$ does not vanish
on $\supp\widetilde{\psi}$ or there is exactly one $x^0\in S^1$ such that \ every $x\in\supp\widetilde{\psi}$
with $\partial_2 f(x)=0$
lies on the curve $\{\delta_r(x^0)\vert\, r>0\}$. 
Now define $\widetilde{f}(x):=f(x)-\partial_1 f(x^0)x_1$. Note
that $\widetilde{f}$ is also $\kappa$-homogeneous of degree one and 
$\nabla\widetilde{f}(x^0)=0$.
Furthermore the analytic family of measures on $S$ will in this case be defined by
$d\sigma_w(x):=G_{\widetilde{f}}(x)^w\psi(x)d\sigma(x)$ (again we take 
$G_{\widetilde{f}}\equiv 1$ if
$\partial_2 f$ does not vanish on $\supp\widetilde{\psi}$). As in the case above, 
the proof of
Theorem \ref{mainthm} is reduced to the estimates \eqref{ft1} and \eqref{ft2} for
 $\alpha>\frac{1}{2}-\frac{1}{h}$.

\section{Auxiliary results on $\kappa$-homogeneous functions}\label{auxiliary}

In this section we shall prove some lemmas which will be useful in the study of the 
oscillatory integrals arizing in the proof of Theorem \ref{mainthm}.

\begin{lemma}\label{localsolutionlemma} Let $U\subseteq\R^2$ be open and
$g\in C^\infty(U)$. If $x^0\in U$ is such that $\partial_2 g(x^0)=0$ and
$\partial_{22}^2 g(x^0)\neq 0$ then there exists a smooth function $\gamma$ of the 
form $\gamma(x_1)=(x_1,\gamma_2(x_1))$, defined in a
neighbourhood of $x_1^0$, such that $\partial_2 g(\gamma(x_1))=0$, and we have

\begin{equation}\label{hessformula}
(g\circ\gamma)''(x_1)=\frac{(\mathrm{Hess\,}g)(\gamma(x_1))}{\partial_{22}^2g(\gamma(x_1))}.
\end{equation}
\end{lemma}

\begin{proof} The existence of $\gamma$ is clear, by the implicit mapping theorem.
But then
\begin{equation*}
(g\circ\gamma)'(x_1)=\partial_1 g(\gamma(x_1)),
\end{equation*}
and hence
\begin{equation*}
(g\circ\gamma)''(x_1)=\partial_{11}^2 g(\gamma(x_1))+
\partial_{21}^2 g(\gamma(x_1))\gamma_2'(x_1).
\end{equation*}
On the other hand, by implicit differentiation we also have the formula
$$\gamma_2'(x_1)=-\frac{\partial_{12}^2 g(\gamma(x_1))}{\partial_{22}^2
 g(\gamma(x_1))}.$$ 
This gives \eqref{hessformula}.
\end{proof}

\begin{lemma}\label{localform}
Let $f\in C^\infty(\R^2\setminus\{0\})$ be
a $\kappa$-homogeneous function. If
$x^0\in\R^2\setminus\{0\}$ with $x^0_1\neq 0$, say $x_1^0>0$, and if
 $n:=\mathrm{ord\,}f(x^0)<\infty,$ then there exists a $\kappa$-homogeneous 
neighbourhood $U$ of $x^0$ on which $f$ is of the form
\begin{equation}\label{taylor}
f(x)=\bigl(x_2-bx_1^{\frac{\kappa_2}{\kappa_1}}\bigr)^ng(x),\quad x\in U,
\end{equation}
where $g$ is a $\kappa$-homogeneous smooth function with $g(x^0)\neq 0$ and $b$
is given by $b:=x_2^0 (x_1^0)^{-\frac{\kappa_2}{\kappa_1}}$. 
If $x^0_2\ne 0,$ then an analogous statement holds, with the roles of the two 
coordinates interchanged.
\end{lemma}

\begin{proof}
By the homogeneity of $f$, we see that $f$ has finite order of 
vanishing  at the point $(1,b)$, too, and that $\partial_2^k f(1,b)\neq 0$, for some 
$k\in\N.$ By  Taylor's formula, applied to the function $y\mapsto f(1,y)$, one therefore gets
$f(1,y)=(y-b)^k G(y)$, with $G(b)\neq 0$. Homogeneity then implies formula 
(\ref{taylor}), with exponent $k$ in place of $n.$ But then evidently 
$k=\mathrm{ord\,}f(x^0)$.
\end{proof}

Recall that a $\kappa$-homogeneous function $f$ of degree $\alpha\in\R$ satisfies 
the following version of Euler's identity:
\begin{equation}\label{euler}
\nabla f(x)\cdot(\kappa_1 x_1,\kappa_2x_2)=\alpha f(x).
\end{equation}

Furthermore it should be noted that the partial derivative $\partial_jf$ of $f$
is also $\kappa$-homogeneous, of degree $\alpha-\kappa_j$.
Applying Euler's identity to the first derivatives of $f$ instead of $f$
then gives

\begin{equation}\label{eqsecder}
D^2 f(x)(\kappa_1 x_1,\kappa_2 x_2)
=((1-\kappa_1)\partial_1f(x),(1-\kappa_2)\partial_2f(x)),
\end{equation}
if $f$ is $\kappa$-homogeneous  of degree $1$.

\begin{lemma}\label{blemma}
Assume that $f$ is $\kappa$-homogeneous of degree $1$ and smooth away from 
the origin, let $x\in\R^2\setminus\{0\}$ and $j\in\{1,2\}$. If $\kappa_j\neq 1$ and 
$\partial_j f(x)\neq 0$, then either $D^2 f(x)$ is non-degenerate, or $\partial^2_{jj}f(x)\neq 0$.
\end{lemma}

\begin{proof}
Assume for instance $j=1$. Since $\partial_1 f(x)\neq 0$ and $\kappa_1\neq 1$, \eqref{eqsecder} implies 
$D^2 f(x)\neq 0$, hence the rank of $D^2 f(x)$ is at least one. Therefore, if $D^2 f(x)$ 
is degenerate, it has rank one. By (\ref{euler}), we have 
\begin{equation*}
\kappa_1 \partial_{11}^2 f(x)x_1+\kappa_2 \partial_{12}^2 f(x)
x_2=(1-\kappa_1)\partial_1f(x)\neq 0.
\end{equation*}

If $\partial_{11}^2f(x)$ were $0$, this would imply $\partial_{12}^2 f(x)\neq 0$, so 
that $D^2 f(x)$ would be non-degenerate. We therefore conclude that 
 $\partial_{11}^2f(x)\neq 0$. The case $j=2$ can be treated in the same way.
\end{proof}

\begin{prop}\label{alemma}
Assume that $f$ is  $\kappa$-homogeneous of
degree $1$, smooth away from the origin, $\mathrm{ord\,}f<\infty$, and that $\kappa_2\neq 1$. 
Let $x^0\in\R^2\setminus\{0\}$.
If $\sigma^0:=-\partial_2 f(x^0)\neq 0$ and $\partial_{22}^2f(x^0)\neq 0$, then
there exists a smooth function  $\gamma$, defined on a neighbourhood $U$ of
$(x_1^0,\sigma^0)$ such that, for all $(x_1,\sigma)\in U$, the function
$x_2\mapsto F(x):=f(x_1,x_2)+\sigma x_2$ has a non-degenerate critical point at
$\gamma(x_1,\sigma)$, and furthermore the function $x_1\mapsto
F(x_1,\gamma(x_1,\sigma))$ does not vanish of infinite order at $x_1^0$.
\end{prop}

\begin{proof} 
By the implicit mapping theorem there exists a smooth function  $\gamma$, defined
on a small neighbourhood of $(x^0_1,\sigma^0)$, such that
\begin{equation*} \partial_2 f(x_1,\gamma(x_1,\sigma))=-\sigma \quad
\mathrm{and}\quad \gamma(x^0_1,\sigma^0)=x^0_2.
\end{equation*}

If we define $\Phi(r,x_1,\sigma):=\delta_r(x_1,\gamma(x_1,\sigma))$,
then homogeneity gives
\begin{equation*} \partial_2 f(\Phi(r,x_1,\sigma))=-r^{1-\kappa_2}\sigma,
\end{equation*}
 and, by differentiating with respect to $x_1$ and $r$ at $r=1$, we obtain
\begin{align*}
\nabla \partial_2 f(x_1,\gamma(x_1,\sigma))\cdot\Phi_{x_1}(1,x_1,\sigma)&=0,\\
\nabla \partial_2 f(x_1,\gamma(x_1,\sigma))\cdot\Phi_r(1,x_1,\sigma)&=
(\kappa_2-1)\sigma\neq 0.
\end{align*}

Since $\Phi_{x_1}(1,x_1,\sigma)=(1,\gamma_{x_1}(x_1,\sigma))\neq 0$, we 
conclude that $\det D_{(r,x_1)}\Phi\neq 0$, and thus, for any fixed $\sigma$, 
the mapping $\Phi$ is a local diffeomorphism w.r.\ to the variables $(r,x_1)$ near 
$r=1, x_1=x_1^0.$

From the homogeneity of the Hessian we also get the identity
\begin{equation*}
(\mathrm{Hess} f)\circ \Phi(r,x_1,\sigma)=r^{2(1-\kappa_1-\kappa_2)}\mathrm{Hess}
f(x_1,\gamma(x_1,\sigma)).
\end{equation*}
The expression on the left hand side cannot vanish of
infinite order with respect to the $(r,x_1)$ variables, due to our assumptions.
 We  therefore see that some derivative of
$\mathrm{Hess} f(x_1,\gamma(x_1,\sigma))$ with respect to $x_1$ does not
vanish at $x_1^0.$ In view of Lemma \ref{localsolutionlemma}, we conclude 
 that the mapping $x_1\mapsto F(x_1,\gamma(x_1,\sigma))$
does not vanish of infinite order at $x_1^0$, because 
$\mathrm{Hess}f=\mathrm{Hess}F$.

\end{proof}

\section{Oscillatory integral estimates}\label{oscestimates}

Throughout this section let $f\in C^\infty(\R^2\setminus\{0\})$ be a real-valued 
$\kappa$-homogeneous function of degree one, and let $\alpha>0$.
By definition (see Section \ref{reduction}),

\begin{equation*}
\widehat{d\sigma_\alpha}(\xi)=e^{-i\xi_3 c}\int_{\R^2} a(x) 
G_f(x)^\alpha e^{-i\xi\cdot(x,f(x))}\,dx,
\end{equation*}
where $a(x)=\psi(x,c+f(x))\sqrt{1+\vert\nabla f(x)\vert^2}$. We assume here that
either $\nabla f$ does not vanish on $\supp a$ or every $x\in\supp a$ with $\nabla f(x)=0$
is of the form $x=\delta_r(x^0)$ with a single $x_0\in S^1$ as indicated in Section \ref{reduction}.
Then $G_f$ is well-defined.
The estimates \eqref{ft1} and \eqref{ft2} actually follow from the next theorem.

\begin{thm}\label{mainprop} Assume that $\kappa\ne (1,1).$ If $h\geq 2$,
$\mathrm{ord}(\mathrm{Hess\,}f)<\infty$ and $\alpha>\frac{1}{2}-\frac{1}{h}$
then, given any bounded neighbourhood $V\subset \R^2$ of the origin, there exist $C,\eps>0$ such that the integral
$J(t,s):=\int_{\R^2}a(x) G_f(x)^\alpha e^{it(f(x)+x\cdot s)}\,dx$ satisfies

\begin{equation}\label{mainpropeq} |J(t,s)|\leq
C\|a\|_{L^1_3}(1+|t|)^{-(\frac{1}{2}+\eps)},\quad\text{for all}\ \ t\in\R,
 s\in\R^2, 
\end{equation}
for any function $a\in\C_0^\infty(V).$
Here, $L^p_k$ denotes the $L^p$-Sobolev space of order $k\in\N.$
\end{thm}

\begin{proof}  We decompose the integral $J(t,s)$ dyadically. To this end,
let $\psi$ be a smooth function supported in the annulus $D=\{1\leq |x|\leq
2\}$ satisfying

\begin{equation*}
\sum_{k=0}^\infty \psi(\delta_{2^k}x)=1,\quad \text{if}\quad 0<|x|<1.
\end{equation*}
After scaling by means of some suitable dilation $\delta_r,$ we may assume  that $V$ is 
 sufficiently small, so that we can write $J(t,s)$ as a
sum of integrals

\begin{equation*} J(t,s)=\sum_{k=0}^\infty J_k(t,s),
\end{equation*} where $J_k$ is defined by
\begin{align*} J_k(t,s)&=\int_{\R^2} a(x)\psi(\delta_{2^k}x) G_f(x)^{\alpha}
e^{it(f(x)+x\cdot s)}\,dx\\ &=2^{-k(\alpha+\kappa_1+\kappa_2)}\int_{\R^2}
a_k(x)\psi(x) G_f(x)^\alpha e^{it2^{-k}(f(x)+x\cdot \sigma)}\,dx,
\end{align*} with $\sigma_1:=2^{k(1-\kappa_1)}s_1$,
$\sigma_2:=2^{k(1-\kappa_2)}s_2$ and $a_k(x):=a(\delta_{2^{-k}}x).$ Using 
Lemma \ref{osclemma} and Theorem \ref{osctheorem} below and
the fact that $\delta:=\alpha+\kappa_1+\kappa_2-\frac{1}{2}-\eps>0,$ if
$\eps>0$ is chosen sufficiently small, we get

\begin{equation*}
\vert J_k(t,s)\vert\leq 2^{-k(\alpha+\kappa_1+\kappa_2)}(1+2^{-k}\vert
t\vert)^{-(\frac{1}{2}+\eps)}\Vert a\Vert_{L^1_3}\leq 2^{-k\delta}(1+\vert
t\vert)^{-(\frac{1}{2}+\eps)}\Vert a\Vert_{L^1_3}
\end{equation*} Thus summation over $k$ yields the desired estimate
\eqref{mainpropeq} of $J(t,s).$
Note that we also used the fact, that $\nabla G_f^\alpha$ is still locally integrable.
\end{proof}

We have thus reduced the estimation of the oscillatory integrals $J(t,s)$ to 
the case where the amplitude function $a$ is supported away from the origin.

\begin{lemma}\label{osclemma}
Let $f$ be $\kappa$-homogeneous of degree one and
$x_0\in\R^2\setminus\{0\}$. For every neighbourhood $U$ of $-\nabla f(x_0)$
and each $N\in\N,$ there exist $C_N>0$ and a compact neighbourhood $K$ of $x_0$
 such that  for all $\sigma\notin U$,
$\lambda\in\R$ and $a\in C_0^\infty(\R^2)$ with $\supp a\subseteq K$

\begin{equation*}
\Bigl\vert\int_{\R^2}a(x)e^{i\lambda(f(x)+\sigma\cdot x)}dx\Bigr\vert\leq
C_N\Vert a\Vert_{L^1_N}(1+\vert\lambda\vert)^{-N}.
\end{equation*}
\end{lemma}

\begin{proof}
The result follows from integrating by parts $N$ times and the trivial estimate
$\vert\int_{\R^2}a(x)e^{i\lambda(f(x)+\sigma\cdot x)}dx\vert\leq\Vert a\Vert_{L^1}$.
\end{proof}

\begin{thm}\label{osctheorem}
Let $f$ be a $\kappa$-homogeneous function of degree one and
$x^0\in\R^2\setminus\{0\}$.
\begin{itemize}
\item[(i)] Fix $j\in\{1,2\}$. If $\kappa_j\neq 1$ and $\partial_jf(x^0)\neq 0$, then 
there are an $\eps>0$ and neighbourhoods $K$ of $x_0$ and $U$ of $\sigma_0:=-\nabla f(x_0)$ such
that, for $\sigma\in U$,
\begin{equation*}
\Bigl\vert\int_{\R^2}a(x)G_f(x)^\alpha e^{i\lambda(f(x)+\sigma\cdot x)}dx\Bigr\vert\leq
C_\alpha\Vert a \Vert_{L^1_3}(1+\vert\lambda\vert)^{-(\frac{1}{2}+\eps)},\quad 
\text{for all}\ \lambda\in\R,
\end{equation*}
for any $\alpha\geq 0$ and any smooth function $a$ supported in $K$.

\item[(ii)] If $\nabla f(x^0)=0,$ and if 
$\alpha>\frac 12-\frac 1h,$ then there are an $\varepsilon>0$ and neighbourhoods
 $K$ of $x^0$ and $U$ of $0$ in $\R^2$ such that, for $\sigma\in U,$

\begin{equation*}
\Bigl\vert\int_{\R^2}a(x) G_f(x)^\alpha
e^{i\lambda(f(x)+\sigma\cdot x)}dx\Bigr\vert\leq C_\alpha\Vert a\Vert_{L^1_3}
(1+\vert\lambda\vert)^{-(\frac{1}{2}+\eps)},\quad 
\text{for all}\ \lambda\in\R,
\end{equation*}
for any smooth function $a$ supported in $K$.

\end{itemize}
\end{thm}

\begin{proof}
In both cases we may assume $\lambda>1$.

(i) Assume $j=2$. By modifying $a$, if necessary, we may assume without loss
of generality that $\alpha=0$ in this case. 
If $D^2 f(x^0)$ is non-degenerate, we apply the stationary phase
method in two dimensions to get the stronger estimate

\begin{equation*}
\Bigl\vert\int_{\R^2}a(x) e^{i\lambda(f(x)+x\cdot\sigma)}dx\Bigr\vert\leq C\Vert
a\Vert_{L^1_2}(1+\vert\lambda\vert)^{-1}.
\end{equation*}
Assume therefore that $D^2f(x^0)$ is degenerate. Then Lemma \ref{blemma} implies $\partial_{22}^2f(x^0)\neq 0$.
That means that we can use the stationary phase method in the $x_2$-variable in order to get

\begin{multline*}
\int_\R a(x)e^{i\lambda(f(x)+\sigma\cdot x)}\,dx_2\\
=\sqrt{2\pi i}\lambda^{-\frac{1}{2}}\frac{a(x_1,\gamma(x_1))}{\sqrt{\partial^2_{22} f(x_1,\gamma(x_1))}}
e^{i\lambda(f(x_1,\gamma(x_1))+\sigma_2
\gamma(x_1)+\sigma_1 x_1)}+R(x_1,\sigma,\lambda),
\end{multline*}
where the remainder term $R$ can be estimated by
$\vert R(x_1,\sigma,\lambda)\vert\leq\lambda^{-1}\Vert
a\Vert_{L^1_2},$ locally uniformly in $\sigma$ and $x_1$. Furthermore, by 
Proposition \ref{alemma}, the phase function 
$x_1\mapsto f(x_1,\gamma(x_1))+\sigma_2\gamma(x_1)+\sigma_1 x_1$ does not vanish
 of infinite order at $x_1^0$. Thus, van der Corput's lemma gives

\begin{equation*}
\Bigl\vert\int_{\R} a(x_1,\gamma(x_1))e^{i\lambda(f(x_1,\gamma(x_1))+\sigma\cdot
(x_1,\gamma(x_1)))}\,dx_1\Bigr\vert\leq C (1+\vert\lambda\vert)^{-\frac{1}{k}}\Vert
a\Vert_{L^1_2},
\end{equation*} with a certain $k\in\N$, which gives the estimate in (i).

(ii) Let us assume without loss of generality that $x^0_1\neq 0$, say  $x^0_1> 0.$
We denote the complete phase function by
$F_\sigma(x):=f(x)+\sigma\cdot x$. By \eqref{euler}, we have  $f(x^0)=0,$ and therefore,
by Lemma \ref{localform},  

\begin{equation*} f(x)=\bigl(x_2-bx_1^{\frac{\kappa_2}{\kappa_1}}\bigr)^n g(x),
\end{equation*}
where $g$ is a smooth function on a $\kappa$-homogeneous neighbourhood of $x^0$ 
with $g(x^0)\neq 0,$ and where $n:=\mathrm{ord}f(x^0)\ge 2$.
This implies that $G_f(x)=\vert x_2-bx_1^{\frac{\kappa_2}{\kappa_1}}\vert^n G(x)$, 
where $G$ is $\kappa$-homogeneous and smooth, with $G(x^0)\neq 0$. Introduce 
the new coordinate

\begin{equation*}
z:=x_2-bx_1^{\frac{\kappa_2}{\kappa_1}}.
\end{equation*}
The map $\Phi(x_1,z):=(x_1,x_2)$ is a local diffeomorphism near $(x_1^0,0)$, since

\begin{equation*}
D\Phi(x_1,z)=\begin{pmatrix} 1&0\\ \ast &1\end{pmatrix}.
\end{equation*}
Expressing $F_\sigma$ in the variables $x_1$ and $z$ gives

\begin{equation*}
(F_\sigma\circ\Phi)(x_1,z)=z^n \widetilde{g}(x_1,z)+z\sigma_2+
b x_1^{\frac{\kappa_2}{\kappa_1}}\sigma_2+\sigma_1 x_1,
\end{equation*}
with $\widetilde{g}:=g\circ\Phi$. Assume for instance that $\sigma_2\ge 0,$ and set 
$\tau:=\sigma_2^{\frac{1}{n-1}}$. We have to estimate 

\begin{equation*}
\int_{\R^2}a(x)G_f(x)^\alpha e^{i\lambda F_\sigma(x)}\,dx
=\int_{\R^2} \widetilde{a}(x_1,z)\vert z\vert^{n\alpha}
e^{i\lambda(\varphi_{x_1,\tau}(z)+b\tau^{n-1}x_1^{\frac{\kappa_2}{\kappa_1}}+
\sigma_1 x_1)}\,d(x_1,z),
\end{equation*}
where $\varphi_{x_1,\tau}(z):=z^n\widetilde{g}(x_1,z)+z\tau^{n-1}$ and 
$\widetilde{a}:=(aG^\alpha)\circ\Phi$.

We examine first the case when $\lambda\tau^n\leq 1,$ in which the theorem can
obviously be reduced to proving
\begin{equation}\label{firstcase_red}
\Bigl\vert\int_{\R}\widetilde{a}(x_1,z)\vert z\vert^{n\alpha}e^{i\lambda\varphi_{x_1,\tau}(z)}dz\Bigr\vert\leq
C\Vert\widetilde{a}\Vert_{L^1_2}\max\bigl(\lambda^{-1},\lambda^{-\frac{1+n\alpha}{n}}\bigr),\quad\lambda\tau^n\leq 1,
\end{equation}
since $\frac 1n+\alpha\geq \frac 1h+\alpha>\frac 12$.
But, notice that we may assume that $|z|$ is small, and thus if $z_c$ is a critical 
point of the phase $\varphi_{x_1,\tau}$ then $z_c\sim\tau.$ Moreover, the $n$-th derivative 
of the phase is bounded from below. Therefore the estimate above is suggested by a 
formal application of van der Corput's lemma. 

In order to give a proof, we split the integral into the integral over the region where 
$|z|\le C\lambda^{-1/n}$ and the integral 
over its complement. The first integral is trivially bounded by a constant times 
$(\lambda^{-1/n})^{n\alpha +1}.$ In order to estimate the second integral, observe that 

\begin{align}\label{estimatesphase}
\vert\varphi'_{x_1,\tau}(z)\vert&\geq C\vert z\vert^{n-1}\quad
\text{if $\vert z\vert\gtrsim\tau$,}
&\vert\varphi''_{x_1,\tau}(z)\vert&\leq C\vert z\vert^{n-2},
\end{align}
since  $\widetilde{g}(x_1^0,0)\neq 0.$
Apparently, if $C$ is chosen sufficiently big, these estimates apply where  
$|z|\ge C\lambda^{-1/n},$ and thus an integration by parts yields 
\begin{eqnarray*}
\Bigl\vert\int_{|z|\ge C\lambda^{-1/n}}\widetilde{a}(x_1,z)\vert z\vert^{n\alpha}
e^{i\lambda\varphi_{x_1,\tau}(z)}dz\Bigr\vert
&\lesssim
\Vert\widetilde{a}\Vert_{L^1_3}\lambda^{-1}\left(\frac{z^{n\alpha}}{z^{n-1}}\Bigr\vert_{C\lambda^{-1/n}}+
\int_{C\lambda^{-1/n}}^1\frac{z^{n\alpha-1}}{z^{n-1}}\, dz\right)\\
&\lesssim \Vert\widetilde{a}\Vert_{L^1_3}\mathrm{max}\bigl(\lambda^{-1},\lambda^{-\frac{1+n\alpha}{n}}\bigr).
\end{eqnarray*}
We remark that, with a little more effort, by iterated integration by parts one can even 
estimate by a constant times $\lambda^{-\frac{1+n\alpha}{n}},$ since 
$\vert\varphi^{(k)}_{x_1,\tau}(z)\vert\leq C_k\vert z\vert^{n-k}$, for $ k\ge 2$, but 
we won't need this.

The main difficulties occur if $\lambda\tau^n>1$. Let $\theta\in C_0^\infty(\R)$ be a cut-off function
with $\theta(z)=1$ if $\vert z\vert\leq r$ and $\supp\theta\subseteq[-2r,2r]$, where $r>0$ 
is chosen such that \eqref{estimatesphase} holds for $|z|\geq r\tau$.

By means of an integration by parts, one finds that \eqref{estimatesphase} implies 
\begin{equation*}
\Bigl\vert\int_{\R}\widetilde{a}(x_1,z)\vert z\vert^{n\alpha}(1-\theta)(z/\tau)
e^{i\lambda\varphi_{x_1,\tau}(z)}dz\Bigr\vert\leq C\lambda^{-\frac{1+n\alpha}{n}}.
\end{equation*}
We are therefore finally left to consider

\begin{multline*}
I(\lambda,\sigma):=\int_{\R^2} \widetilde{a}(x_1,z)\vert z\vert^{n\alpha}\theta(z/\tau)
e^{i\lambda(\varphi_{x_1,\tau}(z)+b\tau^{n-1}x_1^{\frac{\kappa_2}{\kappa_1}}+
\sigma_1 x_1)}\,d(x_1,z)\\
=\tau^{1+n\alpha}\int_{\R^2}a_{\tau}(x_1,z)\vert z\vert^{n\alpha}\theta(z)
e^{i\lambda(\tau^n\widetilde{\varphi}_{x_1,\tau}(z)+b\tau^{n-1} x_1^{\frac{\kappa_2}{\kappa_1}}
+\sigma_1 x_1)}\,d(x_1,z),
\end{multline*}
with $a_\tau\in C^\infty$ given by $a_\tau(x_1,z):=(aG^\alpha)(\Phi(x_1,\tau z))$
and $\widetilde{\varphi}_{x_1,\tau}$ by

\begin{equation*}
\widetilde{\varphi}_{x_1,\tau}(z)=z^n \widetilde{g}(x_1,\tau z)+z.
\end{equation*}

Recall that we may assume that $|\tau|$ is sufficiently small. Then one readily sees that
$\widetilde{\varphi}_{x_1,\tau}$ has only nondegenerate critical points, all bounded
away from the origin, and with Hessians also bounded away from zero.
Consider the local solutions $\gamma$ of the equation
$\widetilde{\varphi}_{x_1,\tau}'(\gamma(x_1,\tau))=0$. Then $|\gamma(x_1,\tau)|\sim 1.$ By the method of stationary phase,
we therefore obtain that

\begin{equation}\label{last}
\int_\R a_\tau(x_1,z)\vert z\vert^{n\alpha}e^{i\lambda\tau^n \widetilde{\varphi}_{x_1,\tau}(z)}\,dz
=\sqrt{2\pi i}\lambda^{-\frac{1}{2}}\tau^{-\frac{n}{2}} e^{i\lambda\tau^n
\widetilde{\varphi}_{x_1,\tau}(\gamma(x_1,\tau))}\widetilde{b}(x_1,\tau)+R(x_1,\tau,\lambda),
\end{equation}
where the remainder term $R$ satisfies an estimate
$\vert R(x_1,\tau,\lambda)\vert\leq C_s(\lambda \tau^n)^{-s},$ for any 
$s<1,$
locally uniformly in $x_1$. Choosing $s=1/2+\varepsilon,$ we see that the contribution 
of the remainder term $R$ to $I(\lambda, \sigma)$ is at most of order 
$\tau^{1+n\alpha}(\lambda\tau^n)^{-\frac 12-\varepsilon}=\lambda^{-\frac 12-\varepsilon}
\tau^{1+n\alpha-\frac n2-\varepsilon n}$. Since $1+n\alpha-\frac n2>0,$ the exponent of 
$\tau$ is positive for $\varepsilon$ sufficiently small and  so we get the right estimate.

In order to estimate the contribution of the main term in \eqref{last}, we shall apply 
the method of stationary phase to the integral of the main term with respect to $x_1.$
 We therefore need to estimate the second derivative of
\begin{equation*}
\psi_{\tau,\sigma_1}(x_1):=\tau^n\widetilde{\varphi}_{x_1,\tau}(\gamma(x_1,\tau))
+b\tau^{n-1}x_1^{\frac{\kappa_2}{\kappa_1}}+\sigma_1 x_1.
\end{equation*}
For convenience $F_\sigma$ will be denoted by $F$ throughout the following.

First we observe, that by definition $\psi_{\tau,\sigma_1}(x_1)=(F\circ\Phi)(x_1,\tau
\gamma(x_1,\tau))$. Setting $\widetilde{F}(x_1,z):=(F\circ\Phi)(x_1,\tau z)$ it is
obvious that $\frac{\partial\widetilde{F}}{\partial z}(x_1,\gamma(x_1,\tau))=\tau^n
\widetilde{\varphi}_{x_1,\tau}'(\gamma(x_1,\tau))=0$ and on the other hand
$\frac{\partial\widetilde{F}}{\partial z}(x_1,\gamma(x_1,\tau))=
\tau\frac{\partial(F\circ\Phi)}{\partial z}(x_1,\tau\gamma(x_1,\tau))$. From that
we conclude

\begin{equation*}
\frac{\partial(F\circ\Phi)}{\partial z}(x_1,\tau\gamma(x_1,\tau))=0,
\end{equation*}
first for every $\tau\neq 0$ and by continuity for $\tau=0$ as well.
Furthermore $\tau^2\frac{\partial^2(F\circ\Phi)}{\partial z^2}(x_1,\tau z)=\frac{\partial^2
\widetilde{F}}{\partial z^2}(x_1,z)=\tau^n\widetilde{\varphi}_{x_1,\tau}''(z)$ which means

\begin{equation*}
\frac{\partial^2(F\circ\Phi)}{\partial z^2}(x_1,\tau z)\neq 0,
\end{equation*}
for $\tau\neq 0$, $z\sim 1$ and $\vert x_1-x_1^0\vert$ small enough. Now Lemma
\ref{localsolutionlemma} yields 

\begin{equation}\label{hessfirst}
\psi''_{\tau,\sigma_1}(x_1)=\frac{(\mathrm{Hess}(F\circ\Phi))(x_1,\tau\gamma(x_1,\tau))}
{\frac{\partial^2(F\circ\Phi)}{\partial z^2}(x_1,\tau\gamma(x_1,\tau))}
\end{equation}
In order to compute the right hand side, we will now express the second partial
derivatives of $F\circ\Phi$ in terms of derivatives of $F$. The chain rule for second derivatives gives

\begin{equation}\label{secder}
(F\circ\Phi)''(x_1,z)=^t\!\![\Phi'(x_1,z)]\cdot F''(\Phi(x_1,z))\cdot\Phi'(x_1,z)+F'(\Phi(x_1,z))\circ D^2\Phi(x_1,z)
\end{equation}
where $D^2\Phi(x_1,z)$ is the second total derivative of $\Phi$ at $(x_1,z)$, understood as a bilinear mapping
on $\R^2\times\R^2$ with values in $\R^2$. If $\Phi_j$ denotes the $j$-th component of $\Phi$, the second
term in \eqref{secder} takes the form

\begin{align*}
F'(\Phi(x_1,z))\circ D^2\Phi(x_1,z)&=\begin{pmatrix}
\partial_1 F\ \partial^2_{11}\Phi_1+\partial_2 F\ \partial^2_{11}\Phi_2 &\partial_1 F\ \partial^2_{12}\Phi_1+\partial_2 F\ \partial^2_{12}\Phi_2\\
\partial_1 F\ \partial^2_{12}\Phi_1+\partial_2 F\ \partial^2_{12}\Phi_2 &\partial_1 F\ \partial^2_{22}\Phi_1+\partial_2 F\ \partial^2_{22}\Phi_2
\end{pmatrix}\\
&=\begin{pmatrix} \partial_2 F(\Phi(x_1,z))\,\partial^2_{11}\Phi_2(x_1,z) &0\\ 0&0\end{pmatrix},
\end{align*}
where the second equality follows from the fact that all second partial derivative of $\Phi$ except for
$\partial^2_{11}\Phi_2$ are vanishing. Using this together with
$\frac{\partial F}{\partial x_2}(\Phi(x_1,\tau\gamma(x_1,\tau)))=\frac{\partial(F\circ\Phi)}
{\partial z}(x_1,\tau\gamma(x_1,\tau))=0$ one observes that

\begin{equation*}
(F\circ\Phi)''(x_1,\tau\gamma(x_1,\tau))=[^t\Phi'\cdot (F''\circ\Phi)\cdot\Phi'](x_1,\tau\gamma(x_1,\tau)),
\end{equation*}
and because $\det\Phi'(x_1,z)=1$ this finally leads to $(\mathrm{Hess}(F\circ\Phi))(x_1,\tau\gamma(x_1,\tau))=
(\mathrm{Hess}F)(\Phi(x_1,\tau\gamma(x_1,\tau)))$. Using \eqref{hessfirst} we conclude that
 
\begin{equation*}
\psi''_{\tau,\sigma_1}(x_1)=\frac{(\mathrm{Hess} F)(\Phi(x_1,\tau\gamma(x_1,\tau)))}
{\frac{\partial^2 F}{\partial x_2^2}(\Phi(x_1,\tau\gamma(x_1,\tau)))}
=\frac{(\mathrm{Hess}f)(\Phi(x_1,\tau\gamma(x_1,\tau)))}
{(\partial^2_{22}f)(\Phi(x_1,\tau\gamma(x_1,\tau)))}.
\end{equation*}
Now $\mathrm{Hess}f$ and $\partial^2_{22}f$ are $\kappa$-homogeneous and thus, by Lemma \ref{localform}, $(\mathrm{Hess}f)\circ\Phi\sim z^k$
and $(\partial^2_{22}f)\circ\Phi\sim z^l$ near $(x_1^0,0)$, where $k=\mathrm{ord\,}(\mathrm{Hess}f)(x^0)$ and
$l=\mathrm{ord\,}(\partial^2_{22}f)(x^0)$. Now $l=n-2\leq k$ and therefore
$\psi''_{\tau,\sigma_1}(x_1)\sim(\tau\gamma(x_1))^m$ for $x_1$ near $x_1^0$, where $m:=k-l\geq 0$.
Moreover the higher derivatives of $\psi_{\tau,\sigma_1}$
are also bounded by a constant times $\tau^m$. We conclude
that the main term in \eqref{last} leads to a contribution of order 
$\tau^{1+n\alpha}\lambda^{-\frac{1}{2}}\tau^{-\frac{n}{2}}
(\lambda \tau^m)^{-\frac 12})^\varepsilon$, for any $0<\varepsilon\leq 1$. Since 
$1+n\alpha-\frac n2>0$, choosing $\varepsilon$ sufficiently small, we get the desired 
estimate of  $I(\lambda,\sigma)$, namely

\begin{equation*}
\vert I(\lambda,\sigma)\vert\lesssim\Vert a\Vert_{L^1_3}\lambda^{-\frac{1+\eps}{2}}.
\end{equation*}

\end{proof}

\section{Sharpness of the condition $p>h$}\label{sharpness}

That $p>h$ is necessary can be seen by means of standard examples. Assume that $c=1$ and that $\psi(x)>1$ if
$\vert x_1\vert,\vert x_2\vert\leq 1$.  We fix $\eps>0$ such that $\vert f(x)\vert\leq 1$ if
$\vert x_1\vert,\vert x_2\vert\leq\eps$. If $\kappa_1,\kappa_2<1,$ we denote for $N>>1$ 
 by $g_N$ the characteristic function of the set
$[-N,N]^2\times[-1,1]\subseteq\R^3$ and define
\begin{equation*}
A_N:=\bigl\{y>1\bigr\vert\,\vert y\vert^{1-\kappa_j}\leq\tfrac{N}{2\eps},\ 
 j=1,2\bigr\}.
\end{equation*}
Let $y\in A_N$ and $\vert x\vert<\frac{N}{2}$. Then for every $x'$ with
$\vert x'_j\vert\leq\eps y^{-\kappa_j}$ we have
$\vert yx'_j\vert\leq\frac{N}{2}$ and $y^{\kappa_j}\vert x'_j\vert\leq\eps$ and therefore
$\vert yf(x')\vert=\vert f(y^{\kappa_1}x'_1,y^{\kappa_2}x'_2)\vert\leq 1$.
Thus $g_N(x-yx',y-y(f(x')+1))=1$, which implies
\begin{equation*}
Mg_N(x,y)\geq
\bigl\vert\{x'\bigr\vert\,\vert x'_j\vert\leq\eps\vert y\vert^{-\kappa_j},\ j=1,2\}\bigr\vert\sim\vert y\vert^{-(\kappa_1+\kappa_2)}.
\end{equation*}
If we now assume $M$ to be bounded on $L^p$ then

\begin{equation*}
\int_{A_N} \vert y\vert^{-(\kappa_1+\kappa_2)p}dy
\end{equation*}
has to stay bounded as $N\to\infty$. But this is of course equivalent to $p>\frac{1}{\kappa_1+\kappa_2}$.

On the other hand if $0\neq x^0\in\supp\psi$ is chosen such that $\mathrm{ord}f(x^0)=n$,
 then  the local expansion
from Lemma \ref{localform} gives $f(x)\sim(x_2-bx_1^{\kappa_2/\kappa_1})^n$ near $x_0,$ 
if for example $x_1^0>0$. That yields $\vert yf(x')\vert\leq 1$ whenever
$\vert x'_2-b{x'_1}^{\kappa_2/\kappa_1}\vert\leq \delta\vert y\vert^{-\frac{1}{n}}$, where 
$\delta>0$ is a sufficiently small constant.
If $\vert x\vert\leq\frac{N}{2}$ and $1<y\leq \delta N$ then, for $\delta$ small enough,
\begin{equation*}
Mg_N(x,y)\geq\bigl\vert\{x'\bigr\vert\,\vert x'_2-b{x'_1}^{\kappa_2/\kappa_1}\vert
\leq \delta y^{-\frac{1}{n}}, \vert x'-x^0\vert\leq \delta\vert x^0\vert\}\bigr\vert\sim y^{-\frac{1}{n}}.
\end{equation*}
Thus the integral $\int_1^{\delta N} y^{-\frac{p}{n}}dy$ has to be bounded
for $N\to\infty$ which implies $p>n$. 
In combination, these estimates show that the condition $p>h$ is necessary, since
 $h=\mathrm{ord}f$, if $\kappa_j\ge 1$ for some $j$.


\def\cprime{$'$} \def\cprime{$'$} \def\cprime{$'$} \def\cprime{$'$}
\providecommand{\bysame}{\leavevmode\hbox to3em{\hrulefill}\thinspace}

\end{document}